\documentclass[12pt, reqno]{amsart}

\usepackage{amssymb}  

\theoremstyle{plain}
\newtheorem{theorem}{Theorem}[section]
\newtheorem{lemma}[theorem]{Lemma}
\newtheorem{proposition}[theorem]{Proposition}
\newtheorem{corollary}[theorem]{Corollary}
\newtheorem{definition}[theorem]{Definition}
\theoremstyle{remark}

\numberwithin{equation}{section}

\newcommand{\seclabel}[1]{\label{sec:#1}} 
\newcommand{\thmlabel}[1]{\label{thm:#1}} 
\newcommand{\lemlabel}[1]{\label{lem:#1}} 
\newcommand{\corlabel}[1]{\label{cor:#1}} 
\newcommand{\prplabel}[1]{\label{prp:#1}} 
\newcommand{\deflabel}[1]{\label{def:#1}} 
\newcommand{\eqnlabel}[1]{\label{eqn:#1}} 

\newcommand{\secref}[1]{\ref{sec:#1}} 
\newcommand{\thmref}[1]{\ref{thm:#1}} 
\newcommand{\lemref}[1]{\ref{lem:#1}} 
\newcommand{\corref}[1]{\ref{cor:#1}} 
\newcommand{\prpref}[1]{\ref{prp:#1}} 
\newcommand{\defref}[1]{\ref{def:#1}} 
\newcommand{\eqnref}[1]{\ref{eqn:#1}} 

\newcommand{\peqref}[1]{(\eqnref{#1})} 

\newcommand{\Aut}{\mathcal{A}ut}		
\newcommand{\Ann}{\mathrm{Ann}}		
\newcommand{\End}{\mathcal{E}nd}	
\newcommand{\ZEnd}{\mathcal{ZE}nd}	

\newcommand{\iv}{^{-1}}				
\newcommand{\Fl}{($F_l$)\ }			
\newcommand{\Fr}{($F_r$)\ }			

\title[F-Quasigroups Isotopic to Groups]
{F-Quasigroups Isotopic to Groups}

\author[T.~Kepka]{Tom\'{a}\v{s}~Kepka$^*$}
\thanks{$^*$Partially supported by the institutional grant
MSM 113200007 and by the Grant Agency of Charles University,
grant \#269/2001/B-MAT/MFF}
\address{Department of Algebra \\
MFF UK, Sokolovsk\'{a} 83 \\
186 75 Praha 8, Czech Republic}
\email{kepka@karlin.mff.cuni.cz}
\author[M.~K.~Kinyon]{Michael~K.~Kinyon}
\address{Department of Mathematical Sciences \\
Indiana University South Bend \\
South Bend, IN 46634 USA}
\email{mkinyon@iusb.edu}
\urladdr{http://mypage.iusb.edu/\symbol{126}mkinyon}
\author[J.~D.~Phillips]{J.~D.~Phillips}
\address{Department of Mathematics \& Computer Science \\
Wabash College \\
Crawfordsville, IN 47933 USA}
\email{phillipj@wabash.edu}
\urladdr{http://www.wabash.edu/depart/math/faculty.html{\#}Phillips}

\date{\today}

\subjclass{20N05}
\keywords{F-quasigroup, Moufang loop, generalized modules}

\begin{document}

\begin{abstract}
In \cite{KKP1} we showed that every loop isotopic to an F-quasigroup is a Moufang loop.
Here we characterize, via two simple identities, the class of F-quasigroups
which are isotopic to groups. We call these quasigroups FG-quasigroups.
We show that FG-quasigroups are
linear over groups. We then use this fact to describe their structure. This gives
us, for instance, a complete description of the simple FG-quasigroups. Finally,
we show an equivalence of equational classes between pointed FG-quasigroups
and central generalized modules over a particular ring.

\end{abstract}

\maketitle

\section{Introduction}
\seclabel{intro}

Let $Q$ be a non-empty set equipped with a binary operation
(denoted multiplicatively throughout the paper). For each
$a \in Q$, the left and right translations $L_a$ and $R_a$ are defined by
$L_a x = ax$ and $R_a x = xa$ for all $x \in Q$. The structure
$(Q,\cdot)$ is called a \emph{quasigroup} if all of the right
and left translations are permutations of $Q$ \cite{Br, Pf}.

In a quasigroup $(Q,\cdot)$, there exist transformations
$\alpha, \beta : Q\to Q$ such that $x \alpha(x) = x = \beta(x) x$ for all $x\in Q$.
A quasigroup $Q$ is called a \emph{left F-quasigroup} if
\[
x \cdot yz = xy \cdot \alpha(x)z  \tag{$F_l$}
\]
for all $x, y, z \in Q$. Dually,
$Q$ is called a \emph{right F-quasigroup}  if 
\[
zy \cdot x = z\beta(x) \cdot yx  \tag{$F_r$}
\]
for all $x,y,z \in Q$.
If $Q$ is both a left F- and right F-quasigroup, then $Q$ is called a
(two-sided) \emph{F-quasigroup} \cite{BF, Go, Ke, KKP1, KKP2, Mu, Sa}.

Recall that for a quasigroup $(Q,\cdot)$ and for fixed $a,b\in Q$, the
structure $(Q,+)$ consisting of the set $Q$ endowed with the
binary operation $+ : Q\times Q\to Q$ defined by $x + y =
R_b\iv x \cdot L_a\iv y$ is called a \emph{principal isotope} of
$(Q, +)$. Here $(Q, +)$ is a quasigroup with neutral element
$0 = ab$, that is, $(Q, +)$ is a \emph{loop} \cite{Br}. 
(Throughout this paper, we will use additive notation for loops,
including groups, even if the operation is not commutative.)

To study any particular class of quasigroups, it is useful to understand
the loops isotopic to the quasigroups in the class. In \cite{KKP1}, we
have shown that every loop isotopic to an F-quasigroup is a Moufang loop.
In this paper, which is in some sense a prequel to \cite{KKP1},
we study the structure of a particular subclass of
F-quasigroups, namely those which are isotopic to groups. 
An F-quasigroup isotopic to a group will be called an
\emph{FG-quasigroup} in the sequel.

A quasigroup $Q$ is called \emph{medial} if $xa \cdot by = xb \cdot ay$
for all $x,y,a,b \in Q$. We see that \Fl and \Fr are generalizations of
the medial identity. The main result of {\S}\secref{basics} is that the
class of FG-quasigroups is axiomatized by two stronger  
generalizations the medial identity. In particular, we will show
(Theorem \thmref{main}) that a quasigroup is an FG-quasigroup if and
only if
\[
xy \cdot \alpha(u) v = x\alpha(u) \cdot yv   \tag{$A$}
\]
and
\[
xy \cdot \beta(u) v = x\beta(u) \cdot yv    \tag{$B$}
\]
hold for all $x,y,u,v$.

In {\S}\secref{FG-linear}, we will show that FG-quasigroups are more
than just isotopic to groups; they are, in fact, linear over groups.
A quasigroup $Q$ is said to be \emph{linear} over a group $(Q,+)$
if there exist $f,g\in \Aut(Q,+)$ and $e\in Q$ such that
$xy = f(x) + e + g(y)$ for all $x,y\in Q$. In {\S}\secref{linear-quasi},
we give necessary and sufficient conditions in terms of $f, g,$ and $e$
for a quasigroup $Q$ linear over a group $(Q,+)$ to be an FG-quasigroup.

In {\S}\secref{structure}, we will use the linearity of FG-quasigroups
to describe their structure. For a quasigroup $Q$, set
$M(Q) = \{a \in Q : xa \cdot yx = xy \cdot ax \; \forall x,y \in Q\}$.
We will show (Proposition \prpref{structure}) that in an FG-quasigroup
$Q$, $M(Q)$ is a medial, normal subquasigroup and $Q/M(Q)$ is a group. In
particular, this gives us a complete description of simple FG-quasigroups
(Corollary \corref{simple}) up to an understanding of simple groups.

In {\S}\secref{forms} we codify the relationship between FG-quasigroups
and groups by introducing the notion of \emph{arithmetic form} for an
FG-quasigroup (Definition \defref{form}). This enables us to show an
equivalence of equational classes between (pointed) FG-quasigroups and
certain types of groups with operators (Theorem \thmref{correspondence}
and Lemma \lemref{homom-forms}). Finally, motivated by this equivalence,
we introduce in {\S}\secref{modules} a notion of \emph{central generalized
module} over an associative ring, and we show an equivalence of 
equational classes between (pointed) FG-quasigroups and central generalized
modules over a particular ring (Theorem \thmref{mod-equiv}). 
In \cite{KKP2}, which is the sequel to \cite{KKP1}, we will examine the more
general situation for arbitrary F-quasigroups and introduce a correspondingly
generalized notion of module.

\section{Characterizations of FG-quasigroups}
\seclabel{basics}

\begin{proposition}
\prplabel{basic}
Let $Q$ be a left F-quasigroup. Then
\begin{enumerate}
\item[1.] \qquad $\alpha\beta = \beta\alpha$ and $\alpha$ is an endomorphism
of $Q$.
\item[2.] \qquad $R_aL_b = L_bR_a$ for $a,b \in Q$ if and only if 
$\alpha(b) = \beta(a)$.
\item[3.] \qquad $R_{\alpha(a)}L_{\beta(a)} = L_{\beta(a)}R_{\alpha(a)}$
for every $a \in Q$.
\end{enumerate}
\end{proposition}

\begin{proof}
For (1): 
$x \cdot \alpha\beta(x)\alpha(x) = \beta(x)x \cdot \alpha\beta(x)\alpha(x)
= \beta(x) \cdot x\alpha(x) = \beta(x)x = x = x\alpha(x)$ and so
$\alpha\beta(x) = \beta\alpha(x)$. Further, $xy \cdot \alpha(x)\alpha(y)
= x \cdot y\alpha(x) = xy = xy \cdot \alpha(xy)$ and $\alpha(x)\alpha(y)
= \alpha(xy)$.

For (2):
If $R_aL_b = L_bR_a$, then $ba = R_aL_b\alpha(b) = L_bR_a\alpha(b) =
b \cdot \alpha(b)a, a = \alpha(b)a$ and $\beta(a) = \alpha(b)$.

Conversely, if $\beta(a) = \alpha(b)$ then $b \cdot xa = bx \cdot \alpha(b)a
= bx \cdot \beta(a)a = bx \cdot a$.

Finally (3), follows from (1) and (2).
\end{proof}

\begin{corollary}
\corlabel{commute}
If $Q$ is an F-quasigroup, then $\alpha$ and $\beta$
are endomorphisms of $Q$, and $\alpha \beta = \beta \alpha$.
\end{corollary}

For a quasigroup $(Q,\cdot)$, if the loop isotope
$(Q,+)$ given by $x + y = L_b\iv x\cdot R_a\iv y$ for all $x,y\in Q$
is a associative (\textit{i.e.}, a group), then 
$L_b\iv x \cdot R_a\iv (L_b\iv y\cdot R_a\iv z)
= L_b\iv (L_b\iv x\cdot R_a\iv y) \cdot R_a\iv z$ for all $x,y,z\in Q$.
Replacing $x$ with $L_b x$ and $z$ with $R_a z$, we have that associativity
of $(Q,\circ)$ is characterized by the equation
\begin{equation}
\eqnlabel{assoc}
x \cdot L_b\iv (R_a\iv y\cdot z) = R_a\iv (x\cdot L_b\iv y)\cdot z
\end{equation}
for all $x,y,z \in Q$, or equivalently,
\begin{equation}
\eqnlabel{assoc2}
L_x L_b\iv R_z R_a\iv = R_z R_a\iv L_x L_b\iv
\end{equation}
for all $x,z in Q$.

\begin{lemma}
\lemlabel{group-iso}
Let $Q$ be a quasigroup. The following are equivalent:
\begin{enumerate}
\item[1.] \qquad Every loop isotopic to $Q$ is a group.
\item[2.] \qquad Some loop isotopic to $Q$ is a group.
\item[3.] \qquad For all $x,y,z,a,b \in Q$, \peqref{assoc} holds.
\item[4.] \qquad There exist $a,b\in Q$ such that \peqref{assoc} holds
for all $x,y,z \in Q$.
\end{enumerate}
\end{lemma}

\begin{proof}
The equivalence of (1) and (2) is well known \cite{Br}. (3) and
(4) simply express (1) and (2), respectively, in the form of
equations.
\end{proof}

\begin{lemma}
\lemlabel{FG-char0}
Let $Q$ be an F-quasigroup. The following are equivalent:
\begin{enumerate}
\item[1.] \qquad $Q$ is an FG-quasigroup,
\item[2.] \qquad $x \beta(a) \cdot (L_b\iv R_a\iv y\cdot z)
= (x \cdot R_a\iv L_b\iv y) \cdot \alpha(b) z$ for all $x,y,z\in Q$.
\end{enumerate}
\end{lemma}

\begin{proof}
Starting with Lemma \lemref{group-iso}, observe that \Fr and \Fl
give $R_a\iv (uv) = R_{\beta(a)}\iv u\cdot R_a\iv v$
and $L_b\iv (uv) = L_b\iv u\cdot L_{\alpha(b)}\iv v$ for all $u,v, \in Q$.
Replace $x$ with $x \beta(a)$ and replace $z$ with $\alpha(b) z$.
The result follows.
\end{proof}

\begin{lemma}
\lemlabel{FG-char}
Let $Q$ be an F-quasigroup and let $a,b \in Q$ be such that $\alpha(b) = \beta(a)$.
Then $Q$ is an FG-quasigroup if and only if
$x \beta(a) \cdot y z = x y \cdot \alpha(b) z$ for all $x,y,z\in Q$.
\end{lemma}

\begin{proof}
By Proposition \prpref{basic}(2), $R_aL_b = L_bR_a$ and so $R_a\iv L_b = L_b R_a\iv$.
The result follows from Lemma \lemref{FG-char0} upon 
replacing $y$ with $R_a L_b y$.
\end{proof}

\begin{proposition}
\prplabel{F-as-FG}
The following conditions are equivalent for an F-quasigroup $Q$:
\begin{enumerate}
\item[1.] \quad $Q$ is an FG-quasigroup,
\item[2.] \quad $x\alpha\beta(w) \cdot yz = xy \cdot \alpha\beta(w)z$ for all $x,y,z,w \in Q$.
\item[3.] \quad There exists $w \in Q$ such that $x\alpha\beta(w) \cdot yz = 
xy \cdot \alpha\beta(w)z$ for all $x,y,z \in Q$.
\end{enumerate}
\end{proposition}

\begin{proof}
For given $w\in Q$, set $a = \alpha(w)$ and $b = \beta(w)$. By Corollary
\corref{commute}, $\alpha(b) = \beta(a)$, and so the result follows from 
Lemma \lemref{FG-char}. 
\end{proof}

The preceding results characterize FG-quasigroups among F-quasigroups.
Thus the F-quasigroup laws together with Proposition \prpref{F-as-FG}(2) 
form an axiom base for FG-quasigroups. Now we turn to the main result
of this section, a two axiom base for FG-quasigroups.

\begin{lemma}
\lemlabel{rearrange}
Let $Q$ be an FG-quasigroup. For all $x,y,u,v\in Q$, $L_x L_y\iv R_v\iv R_u
= R_v\iv R_u L_x L_y\iv$.
\end{lemma}

\begin{proof}
Another expression for \Fr is
$R_v\iv R_u = R_{\beta(u)} R_{R_u\iv v}\iv$, and so the
result follows from \peqref{assoc2}.
\end{proof}

\begin{theorem}
\thmlabel{main}
A quasigroup $Q$ is an FG-quasigroup if and only if the identities
$(A)$ and $(B)$ hold.
\end{theorem}

\begin{proof}
Suppose first that $Q$ is an FG-quasigroup. We first verify the following
special case of ($A$): for all $x,y,u,v\in Q$,
\begin{equation}
\eqnlabel{special}
\alpha(x) y \cdot \alpha(u) v = \alpha(x) \alpha(u) \cdot y v
\end{equation}
Indeed, \Fl implies $y = L_u\iv R_{\alpha(u) v}\iv R_{yv} u$. Using 
this and Lemma \lemref{rearrange}, we compute
\[
\alpha(x) y \cdot \alpha(u) v = R_{\alpha(u) v} L_{\alpha(x)} L_u\iv R_{\alpha(u) v}\iv R_{yv} u
= R_{yv} L_{\alpha(x)} L_u\iv u = \alpha(x) \alpha(u) \cdot yv
\]
as claimed.

Next we verify ($B$). For all $x,y,u,v\in Q$,
\[
\begin{array}{rcll}
x \beta(\alpha(u) y)\cdot (u\cdot vy) &=& x \beta(\alpha(u) y)\cdot (uv\cdot \alpha(u)y) & \text{by \Fl} \\
&=& (x\cdot uv) \cdot \alpha(u)y & \text{by \Fr} \\
&=& (xu\cdot \alpha(x)v) \cdot \alpha(u)y & \text{by \Fl} \\
&=& (xu\cdot \beta(\alpha(u)y))\cdot (\alpha(x)v \cdot \alpha(u)y) & \text{by \Fr} \\
&=& (xu\cdot \beta(\alpha(u)y))\cdot (\alpha(x)\alpha(u) \cdot vy) & \text{by \peqref{special}} \\
&=& xu \cdot (\beta(\alpha(u)y) \cdot vy) & \text{by \Fl}
\end{array}
\]
where we have also used Corollary \corref{commute} in the last step.
Replacing $v$ with $R_y\iv v$ and then $y$ with $L_{\alpha{u}}\iv y$, we have ($B$). 
The proof of ($A$) is similar.

Conversely, suppose $Q$ satisfies ($A$) and ($B$). Obviously, ($A$) implies \Fl
and ($B$) implies ($F_r$), and so we may apply Proposition \prpref{F-as-FG} to get
that $Q$ is an FG-quasigroup.
\end{proof}

\section{Quasigroups linear over groups}
\seclabel{linear-quasi}

Throughout this section, let $Q$ be a quasigroup and $(Q,+)$ a group,
possibly noncommutative, 
but with the same underlying set as $Q$. Assume that $Q$ is linear
over $(Q,+)$, that is, there exist $f,g \in \Aut(Q,+)$, $e \in Q$
such that $xy = f(x) + e + g(y)$ for all $x,y\in Q$.

Let $\Phi\in \Aut(Q,+)$ be given by $\Phi(x) = -e + x + e$ for all $x \in Q$.
If we define a multiplication on $Q$ by
$x\cdot_1 y = f(x) + g(y) + e$ for all $x,y \in Q$, then 
$x\cdot_1 y = f(x) + e -e + g(y) + e = f(x) + e + \Phi g(y)$.
On the other hand, if we define a multiplication on $Q$ by
$x\cdot_2 y = e + f(x) + g(y)$ for all $x,y \in Q$, then
$x\cdot_2 y = \Phi\iv f(x) + e + g(y)$. In particular, there is
nothing special about our convention for quasigroups linear over groups;
we could have used $(Q,\cdot_1)$ or $(Q,\cdot_2)$ instead.

\begin{lemma}
\lemlabel{F-char-linear}
With the notation conventions of this section,
\begin{enumerate}
\item[1.] $Q$ is a left F-quasigroup if and only if $fg = gf$ and 
$-x + f(x) \in Z(Q,+)$ for all $x \in Q$,
\item[2.] $Q$ is a right F-quasigroup if and only if $fg = gf$ and 
$-x + g(x) \in Z(Q,+)$ for all $x \in Q$,
\item[3.] $Q$ is an F-quasigroup if and only if $fg = gf$ and
$-x + f(x), -x + g(x) \in Z(Q,+)$ for all $x \in Q$.
\end{enumerate}
\end{lemma}

\begin{proof}
First, note that $\alpha(u) = -g\iv(e) - g\iv f(u) + g\iv(u)$
and $\beta(u) = f\iv(u) - f\iv g(u) - f\iv(e)$ for all $u\in Q$.

For (1): Fix $u,v,w\in Q$ and set $x = f(u)$ and $y = gf(v)$. We have
\[
u \cdot vw = f(u) + e + gf(v) + g(e) + g^2(w)
\]
and
\[
uv \cdot \alpha(u)w = f^2(u)
+ f(e) + fg(v) + e - gfg^{-1}(e) - gfg^{-1}f(u) + gfg^{-1}(u)
+ g(e) + g^2(w).
\]
Thus \Fl holds if and only if 
\begin{equation}
\eqnlabel{tmp}
x + e + y =
f(x) + f(e) + fgf^{-1}g^{-1}(y) + e
- gfg^{-1}(e) - gfg^{-1}(x) + gfg^{-1}f^{-1}(x)
\end{equation}
for all $x,y \in Q$.

Suppose \Fl holds. Then setting $x = 0$ in \peqref{tmp} yields
$e + y = f(e) + fgf\iv g\iv (y) + e - gfg\iv(e)$
and $x = 0 = y$ yields $-f(e) + e = e - gfg\iv(e)$. Thus
$-f(e) + e + y = fgf\iv g\iv (y) - f(e) + e$ and
$x + e + y = f(x) + e + y - gfg\iv (x) + gfg\iv f\iv (x)$. Setting
$y = -e$ in the latter equality, we get $-f(x) + x = -gfg\iv (x)
+ gfg\iv f\iv (x)$ and hence $-f(x) + x + e + y = e + y - f(x) + x$.
Consequently, $-f(x) + x \in Z(Q,+)$ for all $x \in Q$ and looking again
at the already derived equalities, we conclude that $fg = gf$.

For the converse, suppose $fg = gf$. Then \peqref{tmp}, after some
rearranging, becomes
\[
(-f(x) + x) + e + y = f(e) + y + (e - f(e)) + (- f(x) + x) .
\]
If we also suppose $-x + f(x) \in Z(Q,+)$ for
all $x \in Q$, then the latter equation reduces to a triviality,
and so \Fl holds.

The proof of (2) is
dual to that of (1), and (3) follows from (1) and (2).
\end{proof}

It is straightforward to characterize F-quasigroups among quasigroups
linear over groups for the alternative definitions $(Q,\cdot_1)$ and
$(Q,\cdot_2)$ above. Recalling that $\Phi(x) = e + x - e$, 
observe that if $-z + f(z) \in Z(Q,+)$ for all $z\in Q$, then
$fg = gf$ if and only if $f\Phi g = \Phi gf$. Using this observation
and Lemma \lemref{F-char-linear}(3), we get the following assertion:
$(Q,\cdot_1)$ is an F-quasigroup if and only if $fg = gf$ and
$-x + f(x), -x + \Phi g(x) \in Z(Q,+)$ for all $x \in Q$.
Similarly, $(Q,\cdot_2)$ is an F-quasigroup if and only if $fg = gf$ and
$-x + \Phi^{-1}f(x), -x + g(x) \in Z(Q,+)$ for all $x \in Q$.

\section{FG-quasigroups are linear over groups}
\seclabel{FG-linear}

Let $h$ and $k$ be permutations of a group $(Q,+)$. Define a multiplication on 
$Q$ by $xy = h(x) + k(y)$ for all $x,y \in Q$. Clearly, $Q$ is a quasigroup.

\begin{lemma}
\lemlabel{RF-linear}
Assume that $Q$ is a right F-quasigroup. Then:
\begin{enumerate}
\item[1.] $h(x + y) = h(x) - h(0) + h(y)$ for all $x,y \in Q$.
\item[2.] The transformations $x \mapsto h(x) - h(0)$ and
$x \mapsto -h(0) + h(x)$ are automorphisms of $(Q,+)$.
\end{enumerate}
\end{lemma}

\begin{proof}
We have $\beta(u) = h^{-1}(u - k(u))$ and
$h(h(w) + k(v)) + k(u) = wv \cdot u = w\beta(u) \cdot vu 
= h(h(w) + kh^{-1}(u - k(u))) + k(h(v) + k(u))$ for all $u,v,w \in Q$.
Then $h(x + y) + z = h(x + kh^{-1}(k^{-1}(z) - z)) + k(hk^{-1}(y) + z)$
for all $x,y,z \in Q$. Setting $z = 0$ we get $h(x + y) = h(x + t) 
+ khk^{-1}(y)$ where $t = kh^{-1}k^{-1}(0)$. Consequently,
$h(y) = h(t) + khk^{-1}(y)$ and $khk^{-1}(y) = -h(t) + h(y)$. Similarly,
$h(x) = h(x + t) + khk^{-1}(0) = h(x + t) - h(t) + h(0),
h(x + t) = h(x) - h(0) + h(t)$. Thus, $h(x + y) = 
h(x) - h(0) + h(t) - h(t) + h(y) = h(x) - h(0) + h(y)$.
This establishes (1). (2) follows immediately from (1).
\end{proof}

\begin{lemma}
\lemlabel{LF-linear}
Assume that $Q$ is a left F-quasigroup. Then:
\begin{enumerate}
\item[1.] $k(x + y) = k(x) - k(0) + k(y)$ for all $x,y \in Q$.
\item[2.] The transformations $x \mapsto k(x) - k(0)$ and $x \mapsto -k(0) + k(x)$
are automorphisms of $(Q,+)$.
\end{enumerate}
\end{lemma}

\begin{proof}
Dual to the proof of Lemma \lemref{RF-linear}.
\end{proof}

Now let $Q$ be an FG-quasigroup, $a,b \in Q, h = R_a, k = L_b$ and
$x + y = h^{-1}(x) \cdot k^{-1}(y)$ for all $x,y \in Q$. Then $(Q,+)$ is a group
(every principal loop isotope of $Q$ is of this form), $0 = ba$
and $xy = h(x) + k(y)$ for all $x,y \in Q$. Moreover, by Lemmas \lemref{RF-linear}
and \lemref{LF-linear}, the transformations $f:x \mapsto h(x) - h(0)$ and
$g:x \mapsto -k(0) + k(x)$ are automorphisms of $(Q,+)$. We have
$xy = f(x) + e + g(y)$ for all $x,y \in Q$ where $e = h(0) + k(0) 
= 0 \cdot 0 = ba \cdot ba$.

\begin{corollary}
\corlabel{FG-linear}
Every FG-quasigroup is linear over a group.
\end{corollary}

\section{Structure of FG-quasigroups}
\seclabel{structure}

Throughout this section, let $Q$ be an FG-quasigroup. By Corollary \corref{FG-linear},
$Q$ is linear over a group $(Q,+)$, that is, there exist $f,g\in \Aut(Q,+)$, $e\in Q$
such that $xy = f(x) + e + g(y)$ for all $x,y\in Q$. Recall the definition
\[
M(Q) = \{a \in Q : xa \cdot yx = xy \cdot ax \; \forall x,y \in Q\} .
\]

\begin{lemma}
\lemlabel{M=Z-e}
$M(Q) = Z(Q,+) - e = \{a \in Q: xa \cdot yz = xy \cdot az \; \forall x,y,z \in Q\}$.
\end{lemma}

\begin{proof}
If $a \in M(Q)$, then $f^2(x) + f(e) + fg(a) + e + fg(y) + g(e) + g^2(x)
= xa \cdot yx = xy \cdot ax = f^2(x) + f(e) + fg(y) + e + fg(a) + g^2(x)$ or,
equivalently, $fg(a) + e + z = z + e + fg(a)$ for all $z \in Q$. The latter
equality is equivalent to the fact that $fg(a) + e \in Z(Q,+)$ or
$a \in f\iv g\iv (Z(Q,+) - e) = Z(Q,+) - f\iv g\iv (e) = Z(Q,+) - e$,
since $f\iv g\iv (e) - e \in Z(Q,+)$. We have shown that
$M(Q) \subseteq Z(Q,+) - e$. Proceeding conversely, we show that
$Z(Q,+) - e \subseteq \{a \in Q: xa \cdot yz = xy \cdot az\}$, and the
latter subset is clearly contained in $M(Q)$.
\end{proof}

\begin{corollary}
\corlabel{M-equiv}
The following conditions are equivalent:
\begin{enumerate}
\item[1.] \qquad $M(Q) = Z(Q,+)$.
\item[2.] \qquad $e \in Z(Q,+)$.
\item[3.] \qquad $0 \in M(Q)$.
\end{enumerate}
\end{corollary}

\begin{lemma}
\lemlabel{ab-in-M}
$\alpha(Q) \cup \beta(Q) \subseteq M(Q)$.
\end{lemma}

\begin{proof}
This follows from Theorem \thmref{main}.
\end{proof}

\begin{lemma}
\lemlabel{M-subquasi}
$M(Q)$ is a medial subquasigroup of $Q$.
\end{lemma}

\begin{proof}
If $u,v,w \in Z(Q,+)$ then $(u - e) \cdot (v - e) =
f(u) - f(e) + e + g(v) - g(e) = w - e
\in Z(Q,+) - g(e) = Z(Q,+) - e = M(Q)$. Thus
$M(Q) = Z(Q,+) - e$ (Lemma \lemref{M=Z-e}) is closed
under multiplication, and it is easy
to see that for each $a,b\in Z(Q,+)$, the equations
$(a-e)\cdot (x-e) = b-e$ and $(y-e)\cdot (a-e) = b-e$
have unique solutions $x,y\in Z(Q,+)$.
We conclude that $M(Q)$ is a subquasigroup
of $Q$. Applying Lemma \lemref{M=Z-e} again, $M(Q)$ is medial. 
\end{proof}

\begin{lemma}
\lemlabel{M-normal}
$M(Q)$ is a normal subquasigroup of $Q$, and $Q/M(Q)$ is a group.
\end{lemma}

\begin{proof}
$Z(Q,+)$ is a normal subgroup of the group $(Q,+)$, and if $\rho$ denotes
the (normal) congruence of $(Q,+)$ corresponding to $Z(Q,+)$, it is easy to 
check that $\rho$ is a normal congruence of the quasigroup $Q$, too. Finally,
by Lemma \lemref{ab-in-M}, $Q/M(Q)$ is a loop, and hence it is a group.
\end{proof}

Putting together Lemmas \lemref{M=Z-e}, \lemref{ab-in-M}, \lemref{M-subquasi},
and \lemref{M-normal}, we have the following.

\begin{proposition}
\prplabel{structure}
Let $Q$ be an FG-quasigroup. Then 
$\alpha(Q) \cup \beta(Q) \subseteq M(Q)
= \{a \in Q: xa \cdot yz = xy \cdot az \; \forall x,y,z \in Q\}$, $M(Q)$
is a medial, normal subquasigroup of $Q$, and $Q/M(Q)$ is a group.
\end{proposition}

\begin{corollary}
\corlabel{simple}
A simple FG-quasigroup is medial or is a group.
\end{corollary}

\section{Arithmetic forms of FG-quasigroups}
\seclabel{forms}

\begin{definition}
\deflabel{form}
An ordered five-tuple $(Q,+,f,g,e)$ will be called an
\emph{arithmetic form} of a quasigroup $Q$ if the following conditions
are satisfied:
\begin{enumerate}
\item The binary structures $(Q,+)$ and $Q$ share the same underlying set
(denoted by $Q$ again);
\item $(Q,+)$ is a (possibly noncommutative) group;
\item $f,g \in \Aut(Q,+)$;
\item $fg = gf$;
\item $-x + f(x), -x + g(x)  \in Z(Q,+)$ for all $x \in Q$;
\item $e \in Q$;
\item $xy = f(x) + e + g(y)$ for all $x,y \in Q$.
\end{enumerate}
If, moreover, $e \in Z(Q,+)$, then the arithmetic form will be called
\emph{strong}.
\end{definition}

\begin{theorem}
\thmlabel{FG-forms}
The following conditions are equivalent for a quasigroup $Q$:
\begin{enumerate}
\item[1.] $Q$ is an FG-quasigroup.
\item[2.] $Q$ has at least one strong arithmetic form.
\item[3.] $Q$ has at least one arithmetic form.
\end{enumerate}
\end{theorem}

\begin{proof}
Assume (1). 
From Corollary \corref{FG-linear} and Lemma \lemref{F-char-linear}(3),
we know that for all $a,b \in Q$, $Q$ has an arithmetic form
$(Q,+,f,g,e)$ such that $0 = ba$. Further, by Lemma \lemref{ab-in-M},
$\alpha(Q) \cup \beta(Q) \subseteq M(Q)$. Now, if the elements
$a$ and $b$ are chosen so that $ba \in \alpha(Q) \cup \beta(Q)$
(for instance, choose $a = b = \alpha\beta(c)$ for some $c \in Q$
and use Corollary \corref{commute}), or merely that $ba \in M(Q)$,
then the form is strong by Corollary \corref{M-equiv}. Thus (2) holds.
(2) implies (3) trivially, and (3) implies (1) by
Lemma \lemref{F-char-linear}(3).
\end{proof}

\begin{lemma}
\lemlabel{rigid}
Let $(Q,+,f_1,g_1,e_1)$ and $(Q,*,f_2,g_2,e_2)$ be arithmetic forms of the same 
FG-quasigroup $Q$. If the groups $(Q,+)$ and $(Q,*)$ have the same neutral element
$0$, then $(Q,+) = (Q,*)$, $f_1 = f_2, g_1 = g_2$, and $e_1 = e_2$.
\end{lemma}

\begin{proof}
We have $f_1(x) + e_1 + g_1(y) = xy = f_2(x) * e_2 * g_2(y)$ for all $x,y \in Q$.
Setting $x = 0 = y$, we get $e_1 = e_2 = e$. Setting $x = 0$ we get
$p(y) = e + g_1(y) = e_2 * g_2(y)$ and so $f_1(x) + p(y) = f_2(x) * p(y)$. But
$p$ is a permutation of $Q$ and $p(y) = 0$ yields $f_1 = f_2$. Similarly,
$g_1 = g_2$ and, finally, $(Q,+) = (Q,*)$.
\end{proof}

\begin{theorem}
\thmlabel{correspondence}
Let $Q$ be an FG-quasigroup. Then there exists a biunique correspondence between
arithmetic forms of $Q$ and elements from $Q$. This correspondence
restricts to a biunique correspondence between strong arithmetic forms of
$Q$ and elements from $M(Q)$.
\end{theorem}

\begin{proof}
Combine Corollary \corref{FG-linear}, Lemma \lemref{F-char-linear}(3),
and Corollary \corref{M-equiv}. 
\end{proof}

\begin{lemma}
\lemlabel{homom-forms}
Let $Q$ and $P$ be FG-quasigroups with arithmetic forms 
$(Q,+,f,g,e_1)$ and $(P,+,h,k,e_2)$, respectively. Let $\varphi:Q \to P$ be a mapping such that 
$\varphi(0) = 0$. Then $\varphi$ is a homomorphism of the quasigroups if and
only if $\varphi$ is a homomorphism of the groups, $\varphi f = h\varphi,
\varphi g = k\varphi$ and $\varphi(e_1) = e_2$.
\end{lemma}

\begin{proof}
This generalization of Lemma \lemref{rigid} has a similar proof.
\end{proof}

Denote by $\mathcal{F}_{g,p}$ the equational class (and category) of pointed
FG-quasigroups. That is $\mathcal{F}_{g,p}$ consists of pairs $(Q,a)$, $Q$ being
an FG-quasigroup and $a \in Q$ a fixed element. If $(P,b) \in \mathcal{F}_{g,p}$
then a mapping $\varphi:Q \to P$ is a homomorphism in $\mathcal{F}_{g,p}$
if and only if $\varphi$ is a homorphism of the quasigroups and $\varphi(a) = b$. 
Further, put $\mathcal{F}_{g,m} = \{(Q,a) \in \mathcal{F}_{g,p}: a \in M(Q)\}$.
Clearly $\mathcal{F}_{g,m}$ is an equational subclass (and also a full 
subcategory) of $\mathcal{F}_{g,p}$.

Let $\varphi:Q \to P$ be a homomorphism of FG-quasigroups. For every
$a \in Q$ we have $(Q, \alpha(a)), (P,\alpha\varphi(a)) \in \mathcal{F}_{g,m}$,
and $\varphi\alpha(a) = \alpha\varphi(a)$. Thus $\varphi$ is a homomorphism in
$\mathcal{F}_{g,m}$. Similarly, $(Q,\beta(a)), (P,\beta\varphi(a)) \in \mathcal{F}_{g,m}$
and $\varphi\beta(a) = \beta\varphi(a)$.

Denote by $\mathcal{G}$ the equational class (and category) of algebras 
$Q(+,f,g,f\iv ,g\iv ,e)$ where $(Q,+)$ is a group and conditions (2)-(6) of 
Definition \defref{form} are satisfied. If $P(+,h,k,h\iv ,k\iv ,e_1) \in \mathcal{G}$,
then a mapping $\varphi:Q \to P$ is a homomorphism in $\mathcal{G}$ if and
only if $\varphi$ is a homomorphism of the groups such that 
$\varphi f = h\varphi, \varphi g = k \varphi$ and $\varphi(e) = e_1$.
Finally, denote by $\mathcal{G}_c$ the equational subclass of $\mathcal{G}$
given by $e \in Z(Q,+)$.

It follows from Theorem \thmref{correspondence} and Lemma \lemref{homom-forms}
that the classes $\mathcal{F}_{g,p}$ and $\mathcal{G}$ are equivalent. That
means that there exists a biunique correspondence
$\Phi:\mathcal{F}_{g,p} \to \mathcal{G}$ such that for every
algebra $A \in \mathcal{F}_{g,p}$, the algebras $A$ and $\Phi(A)$ have the same underlying
set, and if $B \in \mathcal{F}_{g,p}$, then a mapping
$\varphi:A \to B$ is an $\mathcal{F}_{g,p}$-homomorphism if and only if
it is a $\mathcal{G}$-homomorphism.

\begin{corollary}
\corlabel{equiv-classes}
The equational classes $\mathcal{F}_{g,p}$ and $\mathcal{G}$ are equivalent.
The equivalence restricts to an equivalence between
$\mathcal{F}_{g,m}$ and $\mathcal{G}_c$.
\end{corollary}

\section{Generalized modules}
\seclabel{modules}

Let $(G,+)$ be a (possibly noncommutative) group. An endomorphism $\varphi\in \End(G,+)$ will
be called \emph{central} if $\varphi(G) \subseteq Z(G,+)$. We denote by
$\ZEnd(G,+)$ the set of central endomorphisms of $(G,+)$. Clearly, the composition 
of central endomorphisms is again a central endomorphism and $\ZEnd(G,+)$ becomes
a multiplicative semigroup under the operation of composition. Furthermore,
if $\varphi \in \ZEnd(G,+)$ and $\psi \in \End(G,+)$ then 
$\varphi + \psi \in \End(G,+)$ where $(\varphi + \psi)(x) = \varphi(x) + \psi(x)$ for all $x \in G$.
Consequently, $\ZEnd(G,+)$ becomes an abelian group under pointwise addition,
and, altogether, $\ZEnd(G(+))$ becomes an associative ring (possibly without unity).

Let $R$ be an associative ring (with or without unity). 
A \emph{central generalized (left)} $R$\emph{-module} will be a group $(G,+)$
equipped with an $R$-scalar multiplication $R \times G \to G$ such that
$a(x + y) = ax + ay, (a + b)x = ax + bx, a(bx) = (ab)x$ and 
$ax \in Z(G,+)$ for all $a,b \in R$ and $x,y \in G$.

If $G$ is a central generalized $R$-module, then define the \emph{annihilator}
of $G$ to be $\Ann(G) = \{a \in R : aG = 0\}$. It is easy to see that $\Ann(G)$ is
an ideal of the ring $R$.

Let $\textbf{S} = \mathbb{Z}[\textbf{x},\textbf{y},\textbf{u},\textbf{v}]$
denote the polynomial ring in four commuting indeterminates
$\textbf{x}, \textbf{y}, \textbf{u}, \textbf{v}$ over the ring $\textbf{Z}$ of integers.
Put $\textbf{R} = S\textbf{x} + S\textbf{y} + S\textbf{u} + S\textbf{v}$.
That is, $\textbf{R}$ is the ideal of $\textbf{S}$ generated by the indeterminates.
On the other hand, $\textbf{R}$ is a commutative and associative ring
(without unity) freely generated by the indeterminates.

Let $\mathcal{M}$ be the equational class (and category) of central generalized
$\textbf{R}$-modules $G$ such that 
$\textbf{x} + \textbf{u} + \textbf{x} \textbf{u} \in \Ann(G)$ and
$\textbf{y} + \textbf{v} + \textbf{y} \textbf{v} \in \Ann(G)$. Further,
let $\mathcal{M}_p$ be the equational class of pointed () objects from $\mathcal{M}$. That is, $\mathcal{M}_p$ consists of ordered pairs $(G,e)$ where $G \in \mathcal{M}$
and $e \in G$.
Let $\mathcal{M}_c$ denote the subclass of centrally pointed objects from
$\mathcal{M}_p$, \textit{i.e.}, $(G,e)\in \mathcal{M}_c$ iff $(G,e)\in \mathcal{M}_p$
and $e\in Z(G,+)$.

\begin{theorem}
\thmlabel{mod-equiv}
The equational classes $\mathcal{F}_{g,p}$ and $\mathcal{M}_p$
are equivalent. This equivalence restricts to an equivalence between
$\mathcal{F}_{g,m}$ and $\mathcal{M}_c$
\end{theorem}

\begin{proof}
Firstly, take $(Q,a) \in \mathcal{F}_{g,p}$. Let $(Q,+,f,g,e)$ be the
arithmetic form of the FG-quasigroup $Q$, such that $a = 0$ in $(Q,+)$.
Define mappings $\varphi, \mu, \psi, \nu : Q\to Q$ by 
$\varphi(x) = -x + f(x)$, $\mu(x) = -x + f^{-1}(x)$,
$\psi(x) = -x + g(x)$ and $\nu(x) = -x + g^{-1}(x)$ for all $x\in Q$.
It is straightforward to check that $\varphi, \mu, \psi, \nu$ are
central endomorphisms of $(Q,+)$, that they commute pairwise, and
that $\varphi(x) + \mu(x) + \varphi\mu(x) = 0$ and 
$\psi(x) + \nu(x) + \psi\nu(x) = 0$ for all $x \in Q$.
Consequently, these endomorphisms generate a commutative subring of the
ring $\ZEnd(Q,+)$, and there exists a (uniquely determined) homomorphism
$\lambda:\textbf{R} \to \ZEnd(Q,+)$ such that 
$\lambda(\textbf{x}) = \varphi$, $\lambda(\textbf{y}) = \psi$, 
$\lambda(\textbf{u}) = \mu$, and $\lambda(\textbf{v}) = \nu$.
The homomorphism $\lambda$ induces an $\textbf{R}$-scalar multiplication on the group
$(Q,+)$ and the resulting central generalized $\textbf{R}$-module will be
denoted by $\bar{Q}$. We have
$\lambda(\textbf{x} + \textbf{u} + \textbf{x} \textbf{u})
= 0 = \lambda(\textbf{y} + \textbf{v} + \textbf{y}\textbf{v})$
and so $\bar{Q} \in \mathcal{M}$. Now define
$\rho : \mathcal{F}_{g,p} \to \mathcal{M}_p$ by 
$\rho(Q,a) = (\bar{Q},e)$, and observe that $(\bar{Q},e) \in \mathcal{M}_c$
if and only if $e \in Z(Q,+)$.

Next, take $(\bar{Q},e) \in \mathcal{M}_p$ and define
$f,g : Q\to Q$ by $f(x) = x + \textbf{x}x$ and $g(x) = x + \textbf{y}x$ for all $x \in Q$.
We have $f(x + y) = x + y + \textbf{x}x + \textbf{x}y
= x + \textbf{x}x + y + \textbf{x}y = f(x) + f(y)$ for all $x,y \in Q$,
and so $f\in \End(Q,+)$. Similarly, $g\in \End(Q,+)$. Moreover, $fg(x) = 
f(x + \textbf{y}x) =
x + \textbf{y}x + \textbf{x}x + \textbf{x}\textbf{y}x =
x + \textbf{x}x + \textbf{y}x + \textbf{y}\textbf{x}x =
gf(x)$, and therefore $fg = gf$. Still further, if we define
$k : Q\to Q$ by $k(x) = x + \textbf{u}x$ for $x\in Q$, then $fk(x) = 
x + (\textbf{x} + \textbf{u} + \textbf{x} \textbf{u})x =
x = kf(x)$, and it follows that $k = f\iv$ and so $f\in \Aut(Q,+)$.
Similarly, $g\in \Aut(Q,+)$. Of course, 
$-x + f(x) = \textbf{x}x \in Z(Q,+)$ and
$-x + g(x) \in Z(Q,+)$. Consequently, $Q$ becomes an FG-quasigroup under the
multiplication $xy = f(x) + e + g(y)$. Define $\sigma : \mathcal{M}_p\to \mathcal{F}_{g,p}$
by $\sigma(\bar{Q},e) = (Q,0)$. Using
Theorem \thmref{correspondence} and Lemma \lemref{homom-forms}, it is easy
to check that the operators $\rho$ and $\sigma$ represent an
equivalence between $\mathcal{F}_{g,p}$ and $\mathcal{M}_p$. Further,
$0 \in M(Q)$ if and only if $e \in Z(Q,+)$, so that the equivalence
restricts to $\mathcal{F}_{g,m}$ and $\mathcal{M}_c$.
\end{proof}

\end{document}